\documentclass[12pt]{article}

\title{Multivariate normal approximation in geometric probability}
\usepackage{amssymb}
\usepackage{cite} 
\usepackage[british]{babel}

\author{Mathew D.~Penrose\footnote{e-mail: \texttt{m.d.penrose@bath.ac.uk}}
 \\
 \normalsize
 Department of Mathematical Sciences,
 University of Bath,\\
\normalsize
 Claverton Down, Bath BA2 7AY, England.
\and Andrew R.~Wade\footnote{e-mail: \texttt{Andrew.Wade@bris.ac.uk}}
\\
\normalsize
 Department of Mathematics,
 University of Bristol,\\
\normalsize
 University Walk, Bristol BS8 1TW, England.}

\date{April 2007}

\newcommand{\eqco}{\setcounter{equation}{0}}
\newcommand{\thco}{\setcounter{theo}{0}}
\newcommand{\prco}{\setcounter{prop}{0}}
\newcommand{\laco}{\setcounter{lemm}{0}}
\newcommand{\deco}{\setcounter{defn}{0}}
\newcommand{\allco}{\eqco  \thco \prco \laco  \deco}

\newcommand{\ud}{\mathrm{d}}
\newtheorem{theo}{Theorem}
\newtheorem{prop}{Proposition}
\newtheorem{lemm}{Lemma}
\newtheorem{defn}{Definition}
\newcommand{\bean}{\begin{eqnarray*}}
\newcommand{\eean}{\end{eqnarray*}}
\newcommand{\bea}{\begin{eqnarray}}
\newcommand{\be}{\begin{equation}}
\newcommand{\ee}{\end{equation}}
\newcommand{\eea}{\end{eqnarray}}
\newcommand{\tod}{\stackrel{{\cal D}}{\longrightarrow}}

\newcommand{\proof}{\noindent \textbf{Proof. }}

\def\la{\lambda}
\def\Exp{{\mathbb{E}}}
\def\bbP{{\mathbb{P}}}
\def\Pr{\mathbb{P}}
\def\Var{{\mathrm{Var}}}
\def\R{\mathbb{R}}
\def\Z{\mathbb{Z}}
\def\N{\mathbb{N}}
\def\1{{\bf 1 }}
\def\li{\ell_\infty}
\def\0{{\bf 0}}
\def\X{{\cal X}}
\def\FF{{\cal F}}
\def\Po{{\cal P}}
\def\L{{\cal L}}
\def\H{{\cal H}}
\def\U{{\cal U}}
\def\W{{\cal W}}

\def\MM{{\cal M}}
\def\NN{{\cal N}}
\def\ZZ{{\cal Z}}
\def\QQ{{\cal Q}}

\def\tT{\tilde T}

\def\card{{\rm  card}}
\def\BB{{\cal B}}
\def\eps{{\varepsilon}}
\def\xil{{\xi_\lambda}}
\def\mulx{{\mu ^\xi_\lambda }}

\begin{document}

\maketitle

\begin{abstract}
Consider a measure $\mu_\lambda = \sum_x \xi_x \delta_x$
where the sum is over points $x$ of a Poisson point process
of intensity $\la$ on a bounded region in $d$-space,
and $\xi_x$ is a functional determined by the
Poisson points near to $x$,
i.e. satisfying an exponential stabilization condition, along
with a moments condition (examples include statistics for proximity graphs,
germ-grain models and random sequential deposition models).  
A known general result
says the $\mu_\la$-measures
(suitably scaled and centred) 
 of disjoint sets in $\R^d$
are asymptotically 
independent normals as $\la \to \infty$; here we give an
$O(\la^{-1/(2d + \eps)})$ bound on the rate of convergence. 
We illustrate our
result with an explicit multivariate central limit theorem for the
nearest-neighbour graph on Poisson points on a finite collection of 
disjoint intervals.
\end{abstract}

\vskip 2mm


\vskip 2mm


\newpage

\section{Introduction}
\allco
There has been considerable recent interest in providing
central limit theorems (CLTs)
for certain functionals in geometric probability
defined on spatial Poisson
point processes.
 Such functionals include those associated
with
  random spatial
graphs such as the minimal-length spanning tree or the nearest-neighbour graph,
 as well as with germ-grain models and random sequential
packing models. These functionals are random variables given by sums of contributions from
points of a Poisson point process in $\R^d$.

A natural extension to random {\em measures} may be 
provided by keeping track of the location of each contribution 
in $\R^d$. In this way one can obtain a random field indexed by test functions
on $\R^d$ or by subsets of $\R^d$. 
For example, one can consider the measure
induced by a Poisson process with a point mass at each Poisson
point equal to  the distance to its nearest-neighbour; then a typical multivariate statistic
induced by this measure is the vector of total edge-lengths of the
nearest-neighbour graph
on Poisson points over a finite collection of disjoint subsets of $\R^d$.

Under certain conditions, it is known \cite{by2,mdpmult,mpconv} that the measures, appropriately
scaled and centred,
of disjoint sets (or of test functions with disjoint supports)
are asymptotically
distributed as indpendent normals in the large-intensity limit.
The object of the present paper is to 
 give bounds on
 rate  of convergence; these bounds are the main contribution
 of the present paper. We illustrate our result with an application 
 to the nearest-neighbour situation mentioned above.

The unifying concept
of {\em stabilization} on Poisson points
has proved
a useful notion
of local dependence in the context of geometric
probability. This says, roughly speaking, that the contribution from a  
Poisson point is unaffected by changes to the configuration
of Poisson points beyond a certain (random) distance.

 The methodology of stabilization 
has been fruitfully employed, in various guises,
to produce univariate CLTs and laws of large numbers 
for random quantities
in many
problems in geometric probability; see e.g.~\cite{by2,KL,penbook,mdpmult,mpconv,py1,py2,py3,py4}.
The techniques used in this context
 include a martingale method
(see for instance \cite{KL}, and \cite{py1}
where the method is presented for general stabilizing functionals
in geometric probability),
the method of moments \cite{by2},
and Stein's method \cite{py4}, which we employ in the present paper.

The multivariate case, in which several collections of
random variables are considered, has also received some attention. Applications
in geometric probability include, for example, the joint normality of certain random spatial graph
functionals
defined over a finite collection of disjoint regions in $\R^d$. 
There are potential applications
to multivariate statistics,
including nonparametric multi-sample
tests (see e.g.~\cite{rr}).

In the present paper, we employ a form
of {\em Stein's method} (see \cite{stein}), which
has the advantage that it can provide
rates of convergence in the CLT. 
In this context, Stein's method 
is a useful tool for establishing
normal approximations and CLTs for
sums of weakly dependent random variables.
In this paper, the weak dependency structure is
provided by the concept of stabilization on Poisson points.

In the univariate case, the method yields normal approximation of the sum of a single collection
of random variables that are `mostly independent', i.e.~exhibiting 
a local dependency structure. This structure may be captured
using dependency graphs.
 This method was first used
in the context of geometric probability by Avram and Bertsimas 
in
\cite{avbert}
(using the normal approximation
error bounds of \cite{baldirinott1989})
 to provide
CLTs for certain random combinatorial structures that
are locally determined in some sense, including
the $j$-th nearest-neighbour graph, and the Delaunay and Voronoi graphs.

Using the sharper normal approximation bounds of \cite{chenshao}, more general
results for univariate normal approximation 
based on Stein's method
for random point measures were given by Penrose and Yukich
in \cite{py4}. That paper is the foundation
for the present work, which is its multivariate analogue.

Multivariate CLTs for random measures in geometric probability have recently been
proved via the method of moments \cite{by2} and also the martingale method
\cite{mdpmult}. In particular, \cite{mdpmult}
also covers lattice processes (such as percolation),
and does not require
`exponential' stabilization, and so admits a larger class of measures.
The advantage of the results in the present paper is that
information on
rates of convergence is provided.

Beyond the context of geometric probability, mulivariate
central limit theory has been well studied. 
Related results include
multivariate central limit theorems
for sums of independent random variables given in \cite{gotze}.
In
\cite{gr,rr}, multivariate normal approximation bounds are given for sums of (locally) dependent
random variables, often chosen in somewhat special ways,
 including certain statistics defined on random graphs. The results in the present 
paper have the advantage
of being 
more generally applicable in geometric probability.

\section{Main result}
\allco

The basic setting follows that of \cite{py4}. Let $d \in \N$. As in
\cite{py4}, we consider marked point processes in $\R^d$ for the sake of
generality. Let $(\MM, \FF_{\MM}, \bbP_{\MM})$ be a probability space (the mark space).
Let $\xi ( x,s;\X)$ be a measurable $[0,\infty)$-valued function defined for all triples
$(x,s; \X)$, where $x \in \R^d$, $s \in \MM$ are such that $(x,s) \in \X$,
where $\X \subset \R^d \times \MM$ is finite. When $(x,s) \in ( \R^d \times \MM) \setminus \X$,
we abbreviate notation and write $\xi ( x,s ; \X)$
instead of $\xi (x,s; \X \cup \{ (x,s) \} )$.

Given $\X \subset \R^d \times \MM$, $a >0$
and $y \in \R^d$, set $y +a\X :=
\{ (y+ax,s) : (x,s) \in \X \}$, i.e.~translation
and scaling act only on the `spatial' part of $\X$.
For all $\lambda >0$ let
\[ \xi_\lambda ( x,s; \X ) :=
\xi (x,s; x + \lambda^{1/d} (-x + \X)). \]
Thus $\xi_\lambda$ is a `scaled-up' version of $\xi$, defined on a scaled-up
version of the (marked) point set $\X$ dilated around $x$.
We say that
$\xi$ is {\em translation invariant} if
$\xi (x+y,s; y+\X)
= \xi (x,s; \X)$ for all $y \in \R^d$, all
$(x,s) \in \R^d \times \MM$ and all
finite $\X \subset \R^d \times \MM$. 
When $\xi$ is translation invariant, the functional $\xil$
simplifies to $\xil(x,s;\X)= \xi(\lambda^{1/d}x,s;\lambda^{1/d} \X)$.

For $q \in [1,\infty]$, let $\| \cdot \|_q$ denote
the $\ell_q$ norm on $\R^d$. In the sequel we will use
 $q=2$ (the Euclidean norm) and $q=\infty$.

Let $\kappa$ be a probability density function on $\R^d$
with compact support $A \subset \R^d$, where $A$ is non-null (i.e.~has non-zero Lebesgue measure). 
We assume throughout that $\kappa$ is bounded with supremum
denoted by $\| \kappa \|_\infty<\infty$.
For all $\lambda>0$
let $\Po_\lambda$ denote a Poisson point process in
$\R^d \times \MM$ with
intensity measure $(\lambda \kappa(x) \ud x) \times \bbP_\MM (\ud s)$. 

We use the following notion of 
exponential stabilization, as given in \cite{py4} (taking the $A_\lambda$ there
to be  $A$
for all $\lambda$). For $x \in \R^d$ and $r>0$, 
let $B_r (x)$ denote the Euclidean ball centred
at $x$ of radius $r$. 
Let $U$ denote a random element
of $\MM$ with distribution $\bbP_\MM$, independent of $\Po_\lambda$.
 
\begin{defn} \label{expstab}
$\xi$ is {\em exponentially stabilizing} with respect to $\kappa$ and $A$
if for all $\lambda \geq 1$ and all
$x \in A$, there exists an almost surely finite random variable
$R:= R(x,\lambda)$, (a {\em radius of stabilization} for $\xi$ at $x$) such that
\bean
\xil ( x,U ; [ \Po_\lambda \cap 
(B_{\lambda^{-1/d}R} (x) \times \MM)] \cup \X ) 
= \xil (x,U; \Po_\lambda \cap (B_{\lambda^{-1/d}R}(x) \times \MM)),
\eean
for all finite $\X \subset (A \setminus B_{\lambda^{-1/d}R}(x)) \times \MM$,
 and moreover 
\bean
\limsup_{t \to \infty} t^{-1} \log \left( \sup_{\lambda \geq 1, x \in A} \Pr [ R(x,\lambda) > t ]
\right) < 0.
\eean
\end{defn}

Roughly speaking, $R(x,\lambda)$ is a radius of stabilization
if the value of $\xi_\lambda$ at $x$ is unaffected by changes
to the configuration of Poisson points outside $B_{\lambda^{-1/d} R }(x)$.
Exponential stabilization is known to hold for many `locally determined' 
functionals
defined on spatial point processes, and in particular in several cases of interest
in geometric probability; see for example \cite{py4}.
Following \cite{py4}, we also make the following
definition.
\begin{defn} \label{momp}
$\xi$ has a moment of order $p>0$ (with respect
to $\kappa$ and $A$) if
\bea
\label{0309b}
\sup_{\lambda \geq 1, x \in A} \Exp \left[ \left| \xi_\lambda ( x,U
; \Po_\lambda ) \right| ^p \right] < \infty. \eea
\end{defn}

For $\lambda >0$, we define the random weighted point measure
$\mulx$ on $\R^d$, induced by $\xi_\lambda$, by
\[ \mu^\xi_\lambda := \sum_{(x,s) \in \Po_\lambda \cap (A \times \MM)} 
\xil ( x,s ; \Po_\lambda ) \delta_x, \]
where $\delta_x$ is the point measure at $x \in \R^d$.

For $\Gamma \subset \R^d$, 
let $\BB (\Gamma)$ denote the set of bounded Borel-measurable
functions on $\Gamma$. For $f \in \BB(\Gamma)$, let $\langle f, \mulx \rangle := \int_{\Gamma} 
f \ud
\mulx$.
Let $\Phi$ denote, as usual,
the standard normal distribution function on $\R$.
We recall the following univariate normal
approximation
result of
Penrose and Yukich (contained in Theorem 2.1 of \cite{py4}).

\begin{prop} \cite{py4}
\label{0309a}
Let $\xi$ be exponentially stabilizing and satisfy the moment condition (\ref{0309b})
for some $p>3$. For $\Gamma$ a non-null Borel subset of $A$, let $f \in \BB (\Gamma)$ and
put $T:= \langle f, \mulx \rangle$.
Then there exists a constant $C\in (0,\infty)$ depending on $d$, $\xi$, $f$, and $\kappa$
such that for all $\lambda \geq 2$,
\bea
\label{11a}
\sup_{t \in \R} \left| \Pr \left[
\frac{T - \Exp[T]}
{(\Var [ T] )^{1/2} }
\leq t
\right]
- \Phi (t) \right|
\leq
C (\log \lambda)^{3d}
\lambda (\Var  [ T ] )^{-3/2} .
\eea
\end{prop}

For fixed $m \in \N$, let $\Gamma_i$, $i=1,\ldots,m$ be non-null
Borel subsets of $A \subset \R^d$.
 For notational simplicity, 
for $i=1,\ldots,m$ and for $f_i \in \BB(\Gamma_i)$ set
$T_i:= \langle f_i, \mulx \rangle = \int_{\Gamma_i} f_i \ud \mulx$. 
These are the quantities of interest to us in the present paper.
By Proposition \ref{0309a}, under appropriate conditions, we have that, individually,
each $T_i$ satisfies a normal approximation result
of the form of (\ref{11a}). 
For the present paper, we will impose one extra 
condition to control variances such as
$\Var[T_i]$.
\begin{itemize}
\item[(A1)] There exist constants $C_i \in (0,\infty)$
such that for each $i$, for all $\lambda$ sufficiently large,
$\Var[T_i] \geq C_i \lambda$.
\end{itemize}
Under assumption (A1), the bound on the rate of convergence on the right of
(\ref{11a}) (in the case $T=T_i$)
becomes $O ( \lambda^{-1/2} (\log \lambda)^{3d} )$ (compare Corollary 2.1
of \cite{py4}),
and in particular (\ref{11a}) yields the central limit theorems
\[ 
\frac{T_i - \Exp[T_i]}
{(\Var [ T_i] )^{1/2} }
 \tod \NN (0,1),\]
 as $\lambda \to \infty$,
 where $\NN(0,1)$ is the standard normal distribution
 on $\R$
 and `$\tod$' denotes convergence in distribution. 
 Condition (A1) is true in many cases.
In Section \ref{vars} we will give some sufficient conditions for (A1) to hold, and discuss
alternative conditions which lead to somewhat stronger versions of (A1). In particular,
it is often possible to show (under appropriate conditions)
that $\lambda^{-1} \Var[T_i] \to \sigma_i^2$ for some
$\sigma_i^2 \in (0,\infty)$, which may be `explicit' (see Section \ref{vars}).

Our main result, Theorem \ref{thm1} below, extends Proposition \ref{0309a}
to give a multivariate central limit theorem for $(T_i :
i=1,\ldots,m)$, centred and scaled, with a bound on the rate of convergence. 
We impose the additional assumptions that
(A1) holds and that the sub-regions $\Gamma_i$ are pairwise
disjoint and satisfy the natural regularity condition (A2) below.
The central difficulty in extending Proposition \ref{0309a}
to a multivariate version is that the $T_i$ are not, in general, independent. However,
with the aid of stabilization we will show that they are `asymptotically independent'
in an appropriate sense.

To state (A2), we introduce some notation.
For measurable $B \subset \R^d$, let $|B|$ denote the ($d$-dimensional) Lebesgue measure of $B$.
Let $\partial B$
denote the boundary of $B$. 
 For $B \subset \R^d$ and $x \in \R^d$
let $d_q(x,B) := \inf_{y \in B} \| x-y\|_q$.
Also, for $B, B' \subset \R^d$ with $B \cap B' = \emptyset$,
let $d_q(B,B') := \inf_{x \in B, y \in B' } \| x-y\|_q$, i.e.~the
 shortest distance (in the $\ell_q$ sense)
 between $B$ and $B'$. For $r>0$,
let $\partial_r (B)$ denote the
$r$-neighbourhood of the boundary of $B \subset \R^d$ in the $\li$ norm, that is the set
$\{ x \in \R^d: d_\infty(x, \partial B )  \leq r \}$.
 
\begin{itemize}
\item[(A2)] For each $i$, $|\partial_r (\Gamma_i) | = O(r)$ as $r \downarrow 0$.
\end{itemize}

Sufficient conditions for (A2) include that each of the $\Gamma_i$ is convex, or
each is the finite union of convex regions (e.g.~polyhedral). We can now state our main result.

\begin{theo} \label{thm1}
Let $\xi$ be exponentially stabilizing and satisfy the moment
condition (\ref{0309b}) for all  $p \geq 1$. 
Let $m \in \N$. Let $\Gamma_1, \Gamma_2,\ldots,\Gamma_m$
be fixed disjoint non-null
Borel subsets of $A$ satisfying (A2). 
For $i=1,\ldots,m$, let
$f_i \in \BB(\Gamma_i)$ and set 
$T_i:= \langle f_i , \mulx \rangle$.
Suppose that (A1) holds.
Let  $\eps>0$.
Then there exists a constant $C\in(0,\infty)$ depending on
$d$, $\xi$, $\kappa$, $\eps$, $\{f_i\}$ and $\{\Gamma_i\}$, such that, for all $\lambda
\geq 1$, 
\bea
\label{0901d}
 \sup_{t_1,\ldots,t_m \in \R}
\left| \Pr \left[ \bigcap_{i=1}^m \left\{ \frac{ T_i-\Exp [T_i]}{( \Var [T_i])^{1/2}}
\leq t_i \right\}  \right] - \prod_{i=1}^m \Phi(t_i)  \right| \leq
C \lambda^{-1/(2d+\eps)}.
\eea
\end{theo}

In particular, from (\ref{0901d}) we obtain the multivariate central limit theorem that says
\bea
\label{ppi}
 \left( \frac{ T_i-\Exp [T_i]}{( \Var [T_i])^{1/2}} : i=1,\ldots,m \right)
\tod \NN ( 0, I_m ),\eea
as $\lambda \to \infty$,
where $\NN(0,I_m)$ is the $m$-dimensional normal distribution
with mean $0$ and covariance matrix given by the identity matrix $I_m$. 
It was already known \cite{mdpmult, mpconv} that under similar conditions to those
of Theorem \ref{thm1} we have (\ref{ppi}), at least when 
$\lambda^{-1} \Var[ T_i]
\to \sigma_i^2$ for some $\sigma_i^2 \in (0,\infty)$. Theorem \ref{thm1}
adds to this by providing a bound on the rate of convergence.

As an example of the application of Theorem \ref{thm1},
one can take 
$f_i= \1_{\Gamma_i}$
for $i=1,2,\ldots,m$, where $\1_{\Gamma}$ is the indicator function of $\Gamma \subset \R^d$.
We indicate some particular applications of Theorem \ref{thm1} in Section \ref{appl}. 
Under additional technical conditions, 
one can say more about the asymptotic behaviour of the variance
terms in (\ref{0901d}); see Section \ref{vars} below.
 \\

\noindent
{\bf Remark.} The relatively slow rate of convergence
in higher dimensions arises
primarily due to the possibility of strongly dependent points 
in the neighbourhood of the interface of adjacent regions. If all
of the $\Gamma_i$ are separated by a strictly positive distance, then our methods 
can be adapted to yield a rate of convergence
of the same order as in the univariate result (Proposition \ref{0309a}), 
that is $O( \lambda^{-1/2} (\log \lambda)^{3d})$. \\  

For ease of presentation, 
we prove Theorem \ref{thm1} in Section \ref{prfs} under the 
conditions that $\xi$ is translation invariant and that
the mark space is degenerate (i.e.~$\MM = \{ 1 \}$), and so from
now on we
suppress any mention of $\MM$. In particular, point sets such as
$\X$ and $\Po_\lambda$ will be treated as 
(their corresponding) subsets of $\R^d$, and we will write $\xil (x;\X)$ rather than
$\xil(x,1;\X)$. The proof
can be adapted for the general marked case, as in \cite{py4}. 

\section{Towards a proof of Theorem \ref{thm1}}
\label{prfs}
\allco

For everything that follows, we assume that $\Gamma_1,\Gamma_2,\ldots,\Gamma_m$ are (arbitrary) non-null
Borel 
subsets of the bounded region $A \subset \R^d$, such that $\Gamma_i \cap \Gamma_j = \emptyset$ for
$i \neq j$, and condition (A2) holds. Also, for each $i$ we have a function $f_i \in \BB(\Gamma_i)$.

For fixed $\alpha>0$, let $s_\lambda := \alpha \lambda^{-1/d} \log \lambda$,
and let
$\Gamma_i^{\rm bd}$ denote the $s_\lambda$ `boundary region' of
$\Gamma_i \subseteq A$, in the sense 
\bea
\label{0914g}
 \Gamma_i^{\rm bd} : = \left\{ x \in \Gamma_i : d_\infty(x, \partial \Gamma_i)
\leq \alpha \lambda^{-1/d} \log {\lambda} \right\}
= \Gamma_i \cap \partial_{s_\lambda} (\Gamma_i) .\eea
The remainder of the set $\Gamma_i$ we
simply call the `interior' and denote by $\Gamma_i^{\rm in}$, where
\bean
 \Gamma_i^{\rm in} := \left\{ x \in \Gamma_i : d_\infty
(x, \partial \Gamma_i) > \alpha \lambda^{-1/d}
  \log {\lambda} \right\} =
\Gamma_i \setminus \partial_{s_\lambda} (\Gamma_i).\eean
As previously mentioned, we assume that $\xi$ is translation invariant, and that $\MM = \{ 1\}$.

Define
\[  T_i^{\rm bd} := \int_{\Gamma_i^{\rm bd}} f_i \ud \mulx;~~\textrm{ and } ~~ T_i^{\rm in}
:= \int_{\Gamma_i^{\rm in}} f_i \ud\mulx , \]
so that $T_i = T_i^{\rm in}
+ T_i^{\rm bd}$.
To prepare for the proof of Theorem \ref{thm1} we need some
auxiliary lemmas.
For the subsequent results, we
will need the following covering of scaled-up Borel
regions $\lambda^{1/d} B \subset \R^d$ 
by cubes
of side $1$. 
 
First we need some more notation.
 Let $\card(\X)$ denote
the cardinality of set $\X$.
For $x \in \R^d$, 
let $Q_x$ denote
the unit-volume $\li$ ball in $\R^d$ with centre $x$ (i.e., the unit
$d$-cube at $x$). 
For a Borel set $B \subseteq A \subset \R^d$, let
\bea
\label{cover1}
\ZZ_\lambda (B) := \left\{ x \in \Z^d: Q_x \cap
\lambda^{1/d} B \neq \emptyset \right\}, 
\eea
 and set $n_\lambda(B) := \card {(\ZZ_\lambda
(B))}$. Then the covering of $\lambda^{1/d} B$
is \bea
\label{cover2}
\QQ_\lambda (B) := 
\{ Q_z  
: z \in \ZZ_\lambda(B) \}.
\eea
 
The next result gives error bounds for approximating the volume of
$\lambda^{1/d} \Gamma_i$ or of $\lambda^{1/d} \Gamma_i^{\rm bd}$
(as defined at (\ref{0914g}))
 by the number of unit cubes in $\Z^d$ in
 its covering (as defined at (\ref{cover1}) and (\ref{cover2})).
 
\begin{lemm} \label{0305a}
Let $\Gamma_i$ be a non-null
Borel subset of $A \subset \R^d$ such that
$| \partial_r (\Gamma_i) | = O(r)$ as $r \downarrow 0$. 
 Then, as $\lambda \to \infty$,
\bea
\label{0901e}
 n_\lambda (\Gamma_i) - |\lambda^{1/d} \Gamma_i| =
O \left(\lambda^{(d-1)/d} \right). 
\eea
 Define $\Gamma_i^{\rm bd}$ as at (\ref{0914g}). Then, 
as $\lambda \to \infty$,
\bea
\label{0901eq}
 n_\lambda (\Gamma_i^{\rm bd}) - |\lambda^{1/d} \Gamma_i^{\rm bd}| =
O \left(\lambda^{(d-1)/d} \log \lambda \right).
\eea
 \end{lemm}
\proof
 There exists a constant $c \in (0,\infty)$ (depending only on $d$) such that,
for any $\lambda>0$, and any  non-null
Borel subset $B$ of $A$,
\bean
 \lambda^{1/d} B \subseteq 
\bigcup_{z \in \ZZ_\lambda(B)}  Q_z
\subseteq \lambda^{1/d}
B \cup \partial_{c} (\lambda^{1/d} B) , \eean
and hence 
\bea
\label{ppp}
 |\lambda^{1/d} B
| \leq n_\lambda (B) \leq |\lambda^{1/d} B| +
 | \partial_{c} (\lambda^{1/d} B)| 
= | \lambda^{1/d} B| 
+  \lambda | \partial_{c \lambda^{-1/d}} (B)|. \eea
In the case $B=\Gamma_i$, the regularity
 assumption 
that $|\partial_r (\Gamma_i)| = O(r)$ 
as $r \downarrow 0$ 
implies that
 $ |\partial_{c \lambda^{-1/d}} (\Gamma_i)| = O( 
  \lambda^{-1/d})$. Thus (\ref{0901e}) follows from (\ref{ppp}). 
 
 In the case $B=\Gamma_i^{\rm bd}$, we have that
 \bean | \partial_{ c \lambda^{-1/d}} (\Gamma_i^{\rm bd}) | 
   \leq  | \partial_{c\lambda^{-1/d}+s_\lambda} (\Gamma_i) | = O (s_\lambda) ,\eean
 as $\lambda \to \infty$, again by the regularity assumption on $\Gamma_i$. Thus
 (\ref{ppp}) yields (\ref{0901eq}) in this case.
$\square$ \\
  
 Once more consider a Borel subset $B$ of $A \subset \R^d$ and the covering $\QQ_\lambda (B)$
of $\lambda^{1/d} B$. 
For all $z \in \ZZ_\lambda(B)$, the number of points of
$\Po_\lambda \cap \lambda^{-1/d} Q_z$ is a Poisson random variable $N_z$ with
parameter   $\nu_z  := \lambda \int_{\lambda^{-1/d} Q_z} \kappa (x)
\ud x$. Assuming $\nu_z>0$, choose an ordering on the points of
$\Po_\lambda \cap \lambda^{-1/d} Q_z$ uniformly at random from all $N_z!$
possibilities. List the points as $X_{z,1}, \ldots, X_{z,N_z}$,
where conditional on the value of $N_z$, the random variables
$X_{z,k}$, $k=1,2,\ldots,N_z$ are i.i.d.~on $\lambda^{-1/d} Q_z$ with a density
$\kappa_z (  \cdot  ) := \kappa (  \cdot  ) / \int_{\lambda^{-1/d} Q_z}
\kappa (x) \ud x$. Thus we have the representation 
\[ \Po_\lambda \cap B = \bigcup_{z \in \ZZ_\lambda(B) }
\bigcup_{k=1}^{N_z}
( \{ X_{z,k}  \} \cap B ) . \]
Then for $f$ in $\BB(B)$, 
we can express $\langle f, \mulx \rangle$
as follows:
\bea
\label{0306a}
\langle f, \mulx \rangle = 
 \sum_{z \in \ZZ_\lambda(B) } \sum_{k=1}^{N_z} \xi_\lambda (
X_{z,k} ; \Po_\lambda ) \cdot f (X_{z,k}) \cdot \1_B (X_{z,k}) .
\eea

For all $z \in \ZZ_\lambda(B)$
and for all $k\in \N$, let
$R_{z,k}$ denote the radius of stabilization
of $\xi$ at $X_{z,k}$
if $1 \leq k \leq N_z$ and let
$R_{z,k}=0$ otherwise. Define the event $E_{z,k} := \{ R_{z,k} \leq \alpha \log \lambda \}$.
We define here the  function $\tilde T (B;f)$ as follows,
the idea being that $\tilde T (B;f)$ is, with high probability,
the same as $\langle f, \mulx \rangle$, but exhibits a much more localized dependency 
structure. Set
\bea
\label{0306b}
\tilde T (B;f) := \sum_{z \in \ZZ_\lambda(B) } \sum_{k=1}^{N_z}
\xi_\lambda (X_{z,k} ; \Po_\lambda ) \cdot
\1_{E_{z,k}}  \cdot f (X_{z,k}) \cdot \1_B (X_{z,k}),
\eea
where we use $\1_E$ to denote the indicator random variable of the event $E$.

Recall that $\Gamma_i$, $i=1,2,\ldots,m$ are disjoint non-null Borel
regions
in $A \subset \R^d$ and
$f_i \in \BB(\Gamma_i)$ for
$i=1,2,\ldots,m$. Then for each $i$,
 $\tT(\Gamma_i;f_i)$
is defined by (\ref{0306b}). 
In the same way as we use the abbreviations $T_i$, 
$T_i^{\rm bd}$ and 
$T_i^{\rm in}$,
we let $\tilde T_i := \tilde T (\Gamma_i;f_i)$, 
$\tT_i^{\rm bd} := \tT(\Gamma_i^{\rm bd};f_i)$, and 
$\tT_i^{\rm in} := \tT(\Gamma_i^{\rm in};f_i)$. Thus
$\tT_i = \tT_i^{\rm bd} + \tT_i^{\rm in}$.

For $z \in \ZZ_\lambda(B)$ let $Y_z(B;f)$ be the contribution
to $\tilde T(B;f)$ from the points in $\lambda^{-1/d} Q_z$, i.e.
\bea
\label{0516q}
Y_z(B;f) := \sum_{k=1}^{N_z} \xi_\lambda ( X_{z,k} ; \Po_\lambda)
\cdot \1_{E_{z,k}} \cdot f (X_{z,k}) \cdot \1_B (X_{z,k}), 
\eea 
so that $\tilde T(B;f) = \sum_{z \in \ZZ_\lambda(B)} Y_z(B;f)$.

Let $A_\lambda$, $\lambda \geq 1$ be a family of Borel subsets
of $A \subset \R^d$.
The next two results show that the moments condition (\ref{0309b})
implies bounds on the moments of  $Y_z(A_\lambda;f)$
for $f \in \BB (A)$. When we come to apply the two lemmas below, 
we will be taking $A_\lambda = \Gamma_i$ 
or $A_\lambda = \Gamma_i^{\rm bd}$.

\begin{lemm} \label{0516t}
Let $A_\lambda$, $\lambda \geq 1$ be a family of Borel subsets
of $A \subset \R^d$.
If (\ref{0309b}) holds for some $p>0$, then there is a constant
$C\in(0,\infty)$ such that for
all $\lambda \geq 1$, all $k \geq 1$ and $z \in \ZZ_\lambda (A_\lambda)$
\bea
\label{0515s}
\Exp [ | \xi_\lambda ( X_{z,k} ;
\Po_\lambda ) \cdot \1_{A_\lambda} (X_{z,k}) \cdot \1_{\{k \leq N_z\}} |^p] \leq C.
\eea
\end{lemm}
\proof
It suffices to consider the case with $A_\lambda=A$ for all $\lambda$.
The proof of the lemma 
closely follows that of Lemma 4.2 in \cite{py4},
although our covering is somewhat different.
In the notation of the proof of Lemma 4.2 of \cite{py4}, we have
$\rho_\lambda=1$ and
   $\nu_i = \nu_{z_i} \equiv \lambda \int_{\lambda^{-1/d} Q_{z_i}} \kappa (x) \ud x
\leq \| \kappa \|_\infty$, where we have written $\ZZ_\lambda(B)=\{z_1,\ldots,z_{n_\lambda(B)}\}$. 
Then, following the argument in \cite{py4}, 
we obtain (\ref{0515s}). $\square$

\begin{lemm} \label{0516u}
Let $A_\lambda$, $\lambda \geq 1$, be a sequence of Borel subsets of $A \subset \R^d$,
and suppose $f \in \BB(A)$.
If (\ref{0309b}) holds for some $p > 1$, then for any $q \in (1,p)$
there 
is a constant
$C\in(0,\infty)$ such that for
all $\lambda \geq 1$ and all $z \in \ZZ_\lambda(A_\lambda)$
\bea
\label{0515w}
\| Y_z(A_\lambda;f) \|_q^q \leq C.
\eea
\end{lemm}
\proof The proof 
closely follows that of Lemma 4.3 in \cite{py4},
again with $\rho_\lambda$ there equal to $1$ (and $\nu_i \leq \| \kappa \|_\infty$). Thus, 
with the use of Lemma \ref{0516t} (and the boundedness of $f$), 
we obtain (\ref{0515w}). $\square$

\begin{lemm} \label{0516a}
Suppose that $\xi$ is exponentially stabilizing and satisfies the moments condition
(\ref{0309b}) for some $p>3$. Then there exists a constant $C\in (0,\infty)$ such that
for all $\lambda \geq 2$
\bea
\label{0914j}
\sup_{t \in \R} \left| \Pr \left[ \frac{
\tT_i  - \Exp [\tT_i ]}{ ( \Var [
\tT_i  ] )^{1/2}} \leq t \right] - \Phi (t) \right|
\leq C \lambda (\Var [ \tilde T_i ] )^{-3/2} (\log \lambda)^{3d}.
\eea
Moreover, (\ref{0914j}) holds with $\tT_i$ replaced by $\tT_i^{\rm in}$ everywhere.
\end{lemm}
\proof The statement for $\tT_i$
 follows from equation (4.18) in \cite{py4} with $\rho_\lambda=O(\log \lambda)$, $q=3$,
and taking  the $A_\lambda$ of \cite{py4} to be $\Gamma_i$. In equation
(4.18) of \cite{py4}, $T'_\lambda$ is the equivalent of
our $\tT_i$, $T_\lambda$ is our $T_i$,
and $S$ is our $(\tT_i -\Exp [\tT_i]) (\Var [\tT_i])^{-1/2}$. The
statement for $\tT_i^{\rm in}$ follows in the same way, this time
taking the $A_\lambda$ of \cite{py4} to be $\Gamma_i^{\rm in}$. $\square$

\begin{lemm} \label{0305b}
Suppose that (\ref{0309b}) holds for some $p >2$. Then 
there exist constants $C_1,C_2,C_3 \in (0,\infty)$ such that, for
all $\lambda \geq 2$, 
\bea
\label{0914n}
 \Var [ \tT_i^{\rm bd} ] & \leq &
C_1 \lambda^{(d-1)/d} ( \log{ \lambda } )^{d+1} 
,
\\ 
\label{0914os}
 \Var [ \tT_i ] & \leq & 
C_2 \lambda ( \log{ \lambda } )^d 
, ~~~\textrm{and}\\
\label{0914o}
 \Var [ \tT_i^{\rm in} ]& \leq &
C_3 \lambda ( \log{ \lambda } )^d 
. \eea
\end{lemm}
\proof
First we prove (\ref{0914n}). Consider the covering
$\QQ_\lambda (\Gamma_i^{\rm bd})$
of $\lambda^{1/d} \Gamma_i^{\rm bd}$ by unit $d$-cubes,
as defined at (\ref{cover2}). 
For $z \in \ZZ_\lambda (\Gamma_i^{\rm bd})$ let $Y_z(\Gamma_i^{\rm bd};f_i)$ be the
contribution to $\tT_i^{\rm bd}$ from the points in
$\lambda^{-1/d} Q_z$, as defined at (\ref{0516q}), that is
\bea
\label{0914p}
 Y_z(\Gamma_i^{\rm bd};f_i) := \sum_{k=1}^{N_z}
\xi_\lambda ( X_{z,k}; \Po_\lambda) 
\cdot \1_{E_{z,k}} 
\cdot f_i (X_{z,k}) \cdot \1_{\Gamma_i^{\rm bd}} (X_{z,k}).
\eea 
 Now, using the
representation $\tT^{\rm bd}_i = \sum_{z \in \ZZ_\lambda (\Gamma_i^{\rm bd})}
Y_z (\Gamma_i^{\rm bd};f_i)$, 
we have
\bea
\label{0914k}
 \Var [
\tilde T_i^{\rm bd} ] = \sum_z \Var [Y_z(\Gamma_i^{\rm bd};f_i)] + \sum_{z \neq
w} \mathrm{Cov} [Y_z(\Gamma_i^{\rm bd};f_i),Y_w(\Gamma_i^{\rm bd};f_i)] .
\eea 
By the assumption that (\ref{0309b}) holds for some $p >2$,
by taking $q=2$ and $A_\lambda = \Gamma_i^{\rm bd}$ in
 Lemma \ref{0516u} we have that
$\Var[Y_z(\Gamma_i^{\rm bd};f_i)] \leq V$, for some constant $V <
\infty$, for all $z \in \ZZ_\lambda(\Gamma_i^{\rm bd})$. So by
the Cauchy-Schwarz inequality 
we have
$\mathrm{Cov}[Y_z(\Gamma_i^{\rm bd};f_i),Y_w(\Gamma_i^{\rm bd};f_i)] \leq V$.
Also, 
$Y_z(\Gamma_i^{\rm bd};f_i)$ and $Y_w(\Gamma_i^{\rm bd};f_i)$ are independent if
$d_2(Q_z ,Q_w ) > 2 \alpha \log \lambda$ (by the definition
of $E_{z,k}$).
Further, given $z$, the number of
$w$ for which
$d_2(Q_z,Q_w) \leq 2 \alpha \log \lambda$ is $O((\log \lambda)^d)$. 
Hence (\ref{0914k}) implies that
\bea
\label{0914l}
\Var [ \tilde T_i^{\rm bd} ] & \leq &   n_\lambda (\Gamma_i^{\rm bd})
( V + O( (\log \lambda)^d)) .
\eea
Then by (\ref{0901eq}) 
we have that
\bea
\label{0914m}
n_\lambda(\Gamma_i^{\rm bd})   
 =  \lambda | \Gamma_i^{\rm bd} | + O \left( \lambda^{(d-1)/d} \log \lambda \right) 
 = O ( \lambda^{(d-1)/d} \log \lambda ), \eea
 using (\ref{0914g}) and (A2).
So from (\ref{0914l})
and (\ref{0914m}) we obtain (\ref{0914n}). 

The proof of (\ref{0914os})
follows similarly, using 
$A_\lambda=\Gamma_i$ for all $\lambda$ in Lemma \ref{0516u} and
(\ref{0901e}) in place of (\ref{0901eq}). 
Finally, (\ref{0914o}) follows from (\ref{0914os}), (\ref{0914n})
and the Cauchy-Schwarz inequality, since $\tilde T_i^{\rm in} = \tilde T_i - \tilde T_i^{\rm bd}$.
$\square$

\begin{lemm}
\label{0308a}
Suppose that $\xi$ is exponentially stabilizing and satisfies the moments condition
(\ref{0309b}) for some $p>3$.
Then  there exists a constant
$C \in (0,\infty)$ 
such that for any $\delta >0$, all $\lambda \geq 2$, 
and any $t \in \R$
\begin{equation}
\label{0318a}
\Pr \left[ \left| \frac{
\tT_i - \Exp[ \tT_i]}{ ( \Var [ \tilde T_i ] )^{1/2}} -t \right| \leq \delta \right] \leq \sqrt{ \frac{2}{\pi}}
\delta
+ C
( \log \lambda )^{3d} \lambda (\Var [ \tilde T_i ])^{-3/2} ,
\end{equation}
and also 
\bea
\label{0916d}
\Pr \left[ \left| \frac{
\tT_i^{\rm in} - \Exp [\tT_i^{\rm in}]}
{ ( \Var [ \tilde T_i ] )^{1/2}} -t \right| \leq \delta \right] 
\leq 2 \sqrt{ \frac{2}{\pi}} \delta 
+ C
( \log \lambda )^{3d} \lambda (\Var [ \tilde T_i ])^{-3/2} 
\nonumber\\
+ \Pr \left[ \left| \frac{\tT_i^{\rm bd} - \Exp [\tT^{\rm bd}_i] }
{( \Var [ \tilde T_i ] )^{1/2}}\right|
> \delta \right]
.
\eea
\end{lemm}
\proof First we prove (\ref{0318a}). For the duration of this proof,
write \[ F(t) = \Pr \left[
\frac{ \tT_i -\Exp [\tT_i]}{( \Var [ \tT_i ])^{1/2}} \leq t
\right] . \] Then we have that for $t \in \R$ and $\delta>0$
\bean
&  & \Pr \left[ \left| \frac{
\tilde T_i - \Exp [\tT_i]}{( \Var [ \tT_i ] )^{1/2}} -t \right| \leq
\delta \right]
  =  F(t+\delta) - F(t-\delta) \\
 & = & \Phi(t+\delta) - \Phi(t-\delta)
+\left[
F(t+\delta)  - \Phi(t+\delta) \right]  - \left[ F(t-\delta) - \Phi(t - \delta)
\right] \\
& \leq & \left| \Phi(t+\delta) - \Phi(t-\delta) \right|
+\left|
F(t+\delta)  - \Phi(t+\delta) \right|  + \left| F(t-\delta) - \Phi(t - \delta)
\right|.
 \eean
Then (\ref{0318a}) follows from the Mean Value Theorem (applied to the first term on the right of the above inequality) and
Lemma \ref{0516a} (applied to the other two terms). Finally,
we have that for $\delta>0$
\[ \Pr \left[ \left| \frac{
\tT_i^{\rm in} - \Exp [\tT_i^{\rm in}]}
{ ( \Var [ \tilde T_i ] )^{1/2}} -t \right| \leq \delta \right] 
\leq 
\Pr \left[ \left| \frac{
\tT_i - \Exp [\tT_i]}
{ ( \Var [ \tilde T_i ] )^{1/2}} -t \right| \leq 2\delta \right] 
+ \Pr \left[ \left| \frac{\tT_i^{\rm bd} - \Exp [\tT^{\rm bd}_i]}
{ ( \Var [ \tilde T_i ] )^{1/2}} \right|
> \delta \right].\]
Then using (\ref{0318a}) yields (\ref{0916d}).
$\square$

\begin{lemm}
\label{0309d} Suppose that the moments condition
(\ref{0309b}) holds for all $p \geq 1$, and condition (A2) holds. 
Let $k$ be an even positive integer. Then 
there exists a constant $C \in (0,\infty)$ (depending on $k$) 
such that for all $\lambda \geq 2$,
\bea
\label{0914q}
 \Exp \left[   \left| \tT_i^{\rm bd}
 -\Exp [\tT_i^{\rm bd}] \right|^k
\right] \leq C \lambda^{k(d-1)/(2d)}
( \log \lambda )^{k(1+d)/2}. 
\eea 
\end{lemm}
\proof
 Again consider the covering 
$\QQ_\lambda (\Gamma_i^{\rm bd})$ of $\lambda^{1/d} \Gamma_i^{\rm bd}$
as defined at (\ref{cover2}).
 For $z  \in \ZZ_\lambda (\Gamma_i^{\rm bd})$,
let $\bar Y_z$ be the contribution to
$\tilde T_i^{\rm bd}-\Exp [\tilde T_i^{\rm bd}]$
from cube $Q_z$, that is $\bar Y_z 
:= Y_z(\Gamma_i^{\rm bd};f_i) - \Exp [Y_z(\Gamma_i^{\rm bd};f_i)]$ where
$Y_z(\Gamma_i^{\rm bd};f_i)$ is given by
(\ref{0914p}).
 Thus, for all $z \in \ZZ_\lambda (\Gamma_i^{\rm bd})$,
$\Exp[ \bar Y_z] = 0$ and $\Var[ \bar Y_z] \leq V$ for constant
$V$, by Lemma \ref{0516u}. 

Let $k$ be an even positive integer. Then
\[ \Exp \left[ \left| \tilde T_i^{\rm bd} -\Exp \tilde T_i^{\rm bd}
 \right|^k
\right] = \sum_{z_1 \in \ZZ_\lambda (\Gamma_i^{\rm bd}) }
\sum_{z_2  \in \ZZ_\lambda (\Gamma_i^{\rm bd}) } \cdots
\sum_{z_k  \in \ZZ_\lambda (\Gamma_i^{\rm bd}) } \Exp \left[ \bar Y_{z_1} \bar
Y_{z_2} \cdots \bar Y_{z_k} \right] . \] The term 
$\Exp \left[ \bar Y_{z_1}
\bar Y_{z_2} \cdots \bar Y_{z_k} \right]$ will vanish if any of
the cubes corresponding to the $\bar Y_{z_j}$ is farther than $2
\alpha \log{\lambda}$ from all the other cubes (since then it will
be independent of the other $\bar Y_{z_j}$ and has expectation
zero). In other words, the term vanishes if the appropriate geometric
graph (in the sense of \cite{penbook})
on $z_1, z_2, \ldots, z_k$ has any
isolated vertices. For a non-zero contribution to the sum, we
require the graph to have no isolated vertices --- so it must have
no more than $k/2$ components. So in effect, there are at most
$k/2$ `free' indices of $(z_1,\ldots,z_k)$. Values that are not `free'
have $O((\log \lambda)^d)$ possible values.

Further, $\Exp \left[ \bar Y_{z_1}
\bar Y_{z_2} \cdots \bar Y_{z_k} \right] \leq C$ for some constant
$C$, by Lemma \ref{0516u} (given the moments condition
(\ref{0309b}) for all $p \geq 1$) and H\"older's inequality. Thus for some other
constant also denoted $C$,
\bean
 \sum_{z_1  \in \ZZ_\lambda (\Gamma_i^{\rm bd}) }
\sum_{z_2  \in \ZZ_\lambda (\Gamma_i^{\rm bd}) } \cdots
\sum_{z_k  \in \ZZ_\lambda (\Gamma_i^{\rm bd}) } \Exp \left[ \bar Y_{z_1} \bar
Y_{z_2} \cdots \bar Y_{z_k} \right] & \leq & C (n_\lambda(\Gamma_i^{\rm
bd}))^{k/2} (\log \lambda)^{kd/2}
\\ & \leq & C \lambda^{k(d-1)/(2d)} (\log \lambda )^{k/2} (\log \lambda)^{kd/2},
\eean
 the final inequality by (\ref{0901eq}), (\ref{0914g}) and (A2). Hence we have (\ref{0914q}). $\square$\\
 
The next lemma says that given condition (A1), we can obtain lower bounds
on the variances of $\tilde T^{\rm in}_i$ and  $\tilde T_i$. 
We will need the following result from \cite{py4} (see
(4.17) therein), which says that if $\xi$ is exponentially stabilizing
and satisfies the moments condition
(\ref{0309b}) for some $p>2$, then
\bea 
\label{0318c} 
\left| \Var [ \tilde T_i ] - \Var
[ T_i] \right| & \leq & C \lambda^{-2}
. \eea 

\begin{lemm}
\label{0309e}
Suppose that  (A1) and (A2) are satisfied, and that 
the moments condition (\ref{0309b}) holds for all $p \geq 1$.
Then there exist constants $C \in (0,\infty)$ and $\lambda_0\in [1,\infty)$ such that 
for all $\lambda \geq \lambda_0$
\bea
\label{0914e}
 \Var [
\tilde T_i ] & \geq & C \lambda, ~~~\textrm{and}
\\
\label{0914f}
\Var [ \tilde T_i^{\rm in} ] & \geq & C \lambda.
\eea
\end{lemm}
\proof
These follow in a straightforward manner from (\ref{0318c}), (A1), 
(\ref{0914n}), (\ref{0914os}) and the Cauchy-Schwarz inequality. $\square$

\begin{lemm}
\label{0308b} 
Suppose that $\xi$ is exponentially stabilizing and satisfies the moments condition
(\ref{0309b}) for all $p \geq 1$. Suppose conditions (A1) and (A2) hold.
Then for any $\varepsilon>0$, there exists  $C \in (0,\infty)$ 
 such that
for all $\lambda \geq 1$
\bea
\label{0916a}
 \left| \Pr \left[ \bigcap_{i=1}^m
 \left\{ \frac{ \tT_i - \Exp [\tT_i]}{ ( \Var [
\tilde T_i ])^{1/2}} \leq t_i \right\} \right] - \prod_{i=1}^m \Phi(t_i)
\right|  
 \nonumber\\
 \leq 
  C \lambda^{-1/(2d+\varepsilon)} \! + \left| \prod_{i=1}^m \Pr \left[
\frac{\tT_i^{\rm in} -\Exp [\tT_i^{\rm in}]} {( \Var
[ \tilde T_i ] )^{1/2}} \leq t_i \right]  - \prod_{i=1}^m \Phi(t_i)
\right| . \eea 
\end{lemm}
\proof
We abbreviate our notation for the duration of the current proof
by setting $\sigma_i := (\Var [ \tT_i])^{1/2}$. Then we have
\bea
\label{0916b}
 &  & \left| \Pr \left[
\bigcap_{i=1}^m \left\{ (\tT_i -\Exp [\tT_i])\sigma_i^{-1}
 \leq t_i \right\}  \right] - \prod_{i=1}^m
\Phi(t_i) \right|  \nonumber\\ & \leq & \left| \Pr \left[ \bigcap_{i=1}^m
\left\{  (\tT_i^{\rm in} -\Exp [\tT_i^{\rm in}]) \sigma_i^{-1}
 \leq t_i \right\} \right]
-\prod_{i=1}^m \Phi(t_i) \right| \nonumber\\  &  &
+ \sum_{i=1}^m
\Pr \left[ (\tT_i^{\rm in} -\Exp [\tT_i^{\rm in}] )\sigma_i^{-1}
 \leq t_i, \; ( \tT_i-\Exp[ \tT_i] )\sigma_i^{-1} > t_i \right]  
\nonumber\\ &  &  + \sum_{i=1}^m
\Pr \left[ (\tT_i^{\rm in} -\Exp[ \tT_i^{\rm in}] )\sigma_i^{-1} >  t_i, \; 
(\tT_i-\Exp [\tT_i])\sigma^{-1} 
\leq t_i \right] 
. \eea
For any $\beta >0$, we have
\bea 
\label{0916c}
&  & \max \left( \Pr \left[ ( \tT_i^{\rm in}
-\Exp [\tT_i^{\rm in}]
)\sigma_i^{-1}  \leq t, \; (\tT_i -\Exp [\tT_i]) \sigma_i^{-1} > t \right]
, \right. \nonumber\\ & & ~~~~ \left.
\Pr \left[ ( \tT_i^{\rm in}
-\Exp [\tT_i^{\rm in}]
)\sigma_i^{-1}  > t, \; (\tT_i -\Exp [\tT_i]) \sigma_i^{-1} \leq t \right]
\right)
\nonumber\\ 
& \leq & \Pr \left[ \left| (
\tT_i^{\rm bd} -\Exp [\tT_i^{\rm bd}]) \sigma_i^{-1}
\right| > \lambda^{-\beta} \right] + \Pr
\left[ \left| ( \tT_i^{\rm in} -\Exp [\tT_i^{\rm in}]) \sigma_i^{-1}
-t \right| \leq \lambda^{-\beta} \right] . \eea
Then, from (\ref{0916b}) and (\ref{0916c}) 
\bea
\label{0515a}
 \left| \Pr \left[ \bigcap_{i=1}^m
\left\{ (\tT_i -\Exp [\tT_i]) \sigma_i^{-1}
 \leq t_i \right\} \right] - \prod_{i=1}^m \Phi(t_i)
\right|    \nonumber\\
\leq  \left| \Pr \left[ \bigcap_{i=1}^m \left\{ 
(\tT_i^{\rm in} -\Exp [\tT_i^{\rm in}]) \sigma_i^{-1}
 \leq t_i \right\} \right]
- \prod_{i=1}^m \Phi(t_i) \right| \nonumber\\ 
  + 2\sum_{i=1}^m \Pr \left[ \left| (  \tT_i^{\rm bd}
-\Exp [\tT_i^{\rm bd}] ) \sigma_i^{-1} \right| >
\lambda^{-\beta} \right] + 2\sum_{i=1}^m \Pr \left[ \left| (
\tT_i^{\rm in} -\Exp [\tT_i^{\rm in}] ) \sigma_i^{-1}
 -t_i \right| \leq \lambda^{-\beta} \right] .
\eea
Since $d_2(\lambda^{1/d} \Gamma_i^{\rm in}, \lambda^{1/d}
\Gamma_j^{\rm in})$ is
 at least $2 \alpha \log {\lambda}$ for $i \neq
j$, $\tilde T_i^{\rm in}$, $1 \leq i \leq m$ is a sequence of
mutually independent random variables, so that
\bea
\label{0914a}
 \Pr \left[ \bigcap_{i=1}^m \left\{
( \tT_i^{\rm in} -\Exp [\tT_i^{\rm in}]) \sigma_i^{-1}
  \leq t_i \right\}  \right] =
 \prod_{i=1}^m \Pr \left[ ( \tT_i^{\rm in} -\Exp [\tT_i^{\rm in}])
\sigma_i^{-1}  \leq t_i \right]  . 
\eea
Also, from Markov's inequality, we have that, for $k \in 2 \N$, 
\bea 
\label{0914c}
\Pr \left[ \left|
( \tT_i^{\rm bd} -\Exp [\tT_i^{\rm bd}])
\sigma_i^{-1} \right| > \lambda^{-\beta} \right]
\leq  \Exp \left[
\left| \tT_i^{\rm bd} -\Exp [\tT_i^{\rm bd}]
 \right|^k \right]  \left( \Var [
\tilde T_i ] \right)^{-k/2} \lambda^{k \beta} .
\eea
Then we obtain, from (\ref{0914c}),
with (\ref{0914q}) and
(\ref{0914e}),
\bea
\label{0914b}
\Pr \left[ \left| (
\tT_i^{\rm bd} -\Exp [\tT_i^{\rm bd}])
\sigma_i^{-1} \right|
> \lambda^{-\beta} \right] \leq C \lambda^{k ( \beta -
1/(2d) )} (\log \lambda)^{k(1+d)/2};
\eea
this then gives a bound for the penultimate sum in (\ref{0515a}). To bound the
final sum in (\ref{0515a}), taking $\delta = \lambda^{-\beta}$
we have from (\ref{0916d}), (\ref{0914b}) and (\ref{0914e}) that
\bea
\label{0914bb}
\Pr \left[ \left| (\tT_i^{\rm in} -\Exp [ \tT_i^{\rm in} ]) \sigma_i^{-1}
-t_i \right| \leq \lambda^{-\beta} \right] \nonumber\\
\leq 2 \sqrt{ \frac{2}{\pi}} \lambda^{-\beta}
+ C (\log \lambda)^{3d} \lambda^{-1/2} + C \lambda^{k(\beta-1/(2d))} (\log \lambda)^{k(1+d)/2}.\eea
To obtain the best rates of convergence via this method,
we want to maximize the lowest power of $\lambda^{-1}$
on the right-hand sides of (\ref{0914b}) and (\ref{0914bb}). So we choose $\beta$ such that $-\beta = k\left( \beta -
1/(2d)\right)$, that is, take 
\bea
\label{0914d} 
\beta = \frac{k}{2d(k+1)} .
\eea
For any $\eps >0$
we can choose $k$ large enough in (\ref{0914d})
to give $1/(2d) > \beta \geq
1/(2d+\varepsilon/2)$. Then, for
$\lambda$ sufficiently large, 
$\lambda^{-1/(2d +\varepsilon)} \geq \lambda^{-1/(2d +\varepsilon/2)} (\log \lambda)^{k(1+d)/2}$. 
Now from (\ref{0515a}) and (\ref{0914a}), with the
bounds (\ref{0914b}) and (\ref{0914bb})
  we obtain (\ref{0916a}). 
This completes
the proof of the lemma. $\square$ \\

\noindent
{\bf Proof of Theorem \ref{thm1}.}
To complete the proof we proceed in a similar manner to \cite{py4}. Let
\[ 
E_{\lambda} := \bigcap_{i=1}^m \bigcap_{z \in \ZZ_\lambda(\Gamma_i)} \bigcap_{k=1}^{N_z} 
E_{z,k} , \]
recalling the definition of the event $E_{z,k}$   
just below (\ref{0306a}). By standard Palm theory (e.g.~Theorem 1.6 in
\cite{penbook}) and exponential stabilization (see (4.11) in
\cite{py4}), we have that $\Pr[
E_\lambda^c ] \leq C \lambda ^{-3}$ for $\alpha$ sufficiently
large and some $C\in(0,\infty)$. 
Then $| \tilde T_i -T_i|=0$ except possibly on the set
$E_\lambda^c$, which has probability less than $C \lambda^{-3}$.

For
$i=1,\ldots,m$, let $K_i := (\Var [\tilde T_i])^{-1/2}(\tilde T_i - \Exp
[\tilde T_i])$ and $Z_i := (\Var [\tT_i] )^{-1/2} (T_i-\Exp[T_i])$. 
Then for $\delta>0$  
we have that for any $t_i \in \R$
\[ \{ (Z_i \leq t_i) \Delta (K_i \leq t_i ) \} \subseteq
\{ | K_i -t_i | \leq \delta \} \cup \{ | Z_i -K_i | \geq \delta \},
\]
so that 
\bea \label{0516v} 
& & \left| \Pr \left[
 \bigcap_{i=1}^m \left\{ Z_i \leq t_i \right\} 
\right] - \prod_{i=1}^m \Phi(t_i)
 \right| \leq \left| \Pr 
\left[ \bigcap_{i=1}^m \left\{ K_i \leq t_i \right\} \right] 
- \prod_{i=1}^m \Phi(t_i) \right|  \nonumber\\ &  & \qquad \qquad
 +\sum_{i=1}^m \Pr \left[ \left| K_i -t_i \right| \leq \delta \right]
+ \sum_{i=1}^m \Pr \left[
\left| Z_i-K_i \right| \geq \delta \right] 
.\eea
Then, using (\ref{0916a}) for the first term on the 
right-hand side of the inequality in
 (\ref{0516v}), and (\ref{0318a}) with (\ref{0914e}) for the second, we obtain
\bea
\label{0321a}
 \left| \Pr \left[
 \bigcap_{i=1}^m \left\{ Z_i \leq t_i \right\} 
\right] - \prod_{i=1}^m \Phi(t_i)
 \right| 
 \leq   C \lambda^{-1/(2d+\eps)}  + 
\left| \prod_{i=1}^m
\Pr \left[ \frac{ \tT^{\rm in}_i
- \Exp[
\tT_i^{\rm in}]}
{(\Var[\tT_i])^{1/2}}
\leq t_i \right] - \prod_{i=1}^m \Phi(t_i)
\right| \nonumber\\
+ C \delta    + C (\log \lambda)^{3d} \lambda^{-1/2}
+ \sum_{i=1}^m \Pr \left[ \left| Z_i-K_i
\right| \geq
\delta \right] . \eea
We now consider the second term on the right-hand side of (\ref{0321a}).
For ease of notation, write \[
G_i(t) := \Pr \left[ \frac{ \tilde T_i^{\rm in} -\Exp [\tilde T_i^{\rm in}]}
{(\Var [\tT_i])^{1/2}} \leq t \right], \]
for $i=1,\ldots,m$.
Then \bea \label{0324a}
 \left| \prod_{i=1}^m G_i(t_i) 
- \prod_{i=1}^m \Phi(t_i) \right| 
& = & \left| \sum_{i=1}^m \left[ G_i(t_i) - \Phi(t_i) \right] \prod
_{j=i+1}^m \Phi(t_j) \prod_{k=1}^{i-1} G_k(t_k) \right| \nonumber\\
& \leq & \sum_{i=1}^m \left| G_i(t_i) -\Phi(t_i) \right| 
 . \eea
Writing
\[ H_i(t) := \Pr \left[ \frac{ \tT_i^{\rm in}-\Exp [\tT_i^{\rm in}]}
{( \Var [ \tilde T_i^{\rm in} ])^{1/2}} \leq t \right] ,\]
we have that, for $i=1,\ldots,m$ \bea
\label{0916e}
 \left| G_i(t_i) - \Phi(t_i) \right| \leq \left| H_i ( t_i (1+\gamma_i) )
- \Phi(t_i(1+\gamma_i)) \right| + \left| \Phi(t_i(1+\gamma_i)) -
\Phi(t_i) \right|, \eea
where $1+\gamma_i := \left( \frac{ \Var[ \tilde T_i]} { \Var[\tilde T_i^{\rm in}]}
\right)^{1/2}$. Then,
using Lemma \ref{0516a}
we have
 that the first term on the right-hand side of (\ref{0916e})
satisfies
 \bea
\label{0517b}
\left| H_i(t_i(1+\gamma_i)) - \Phi(t_i(1+\gamma_i)) \right|
 \leq  C (\log \lambda)^{3d} \lambda (\Var [ \tilde T_i^{\rm in}])^{-3/2}
 \leq C \lambda^{-1/2} (\log \lambda)^{3d} , \eea
 by (\ref{0914f}).
In order to deal with the second term on the right-hand side
of (\ref{0916e}), we
need to estimate $\gamma_i$.
We note that
\[ \frac{ \Var [\tilde T_i ]}{ \Var [ \tilde T_i^{\rm in}]}
= 1 + \frac{ \Var [\tilde T_i^{\rm bd} ]}{ \Var [ \tilde T_i^{\rm in}]}
+\frac{2 {\rm Cov} [ \tilde T_i^{\rm in}, \tilde T_i^{\rm bd} ]}{\Var
[\tilde T_i^{\rm in}]} .\]
Then using the upper and lower
variance bounds (\ref{0914n}), (\ref{0914o}), 
(\ref{0914f}), and
the Cauchy-Schwarz
inequality, 
yields
\[ \frac{ \Var [\tT_i ]}{ \Var [ \tT_i^{\rm in}]}
= 1 + O( \lambda^{-1/(2d)} (\log \lambda)^{(2d+1)/2}), \]
so that 
\bea
\label{0517a}
 \gamma_i = O( \lambda^{-1/(2d)} (\log \lambda)^{(2d+1)/2}).
\eea
Since for all $s \leq t$ we have $| \Phi(s) -\Phi(t)|
\leq (t-s) \sup_{s \leq u \leq t} \varphi(u)$
(where $\varphi$ is the standard normal density function), 
we have
\bea
\label{0719x}
&  & \sup_{t_i} | \Phi(t_i(1+\gamma_i))-\Phi(t_i)| \nonumber\\
& \leq & C \sup_{t_i}
\left( |t_i| \lambda^{-1/(2d)} (\log \lambda)^{(2d+1)/2}
\sup_{|u-t_i| \leq t_i C \lambda^{-1/(2d)} (\log \lambda)^{(2d+1)/2}}
\varphi(u) \right) \nonumber\\
& \leq & C \lambda^{-1/(2d)} (\log \lambda)^{(2d+1)/2}.
\eea
So, for the second term on the right-hand side in (\ref{0321a}),
we obtain from (\ref{0324a}), (\ref{0916e}),
(\ref{0517b}) and (\ref{0719x})
\bea
\label{0917b}
& & \sup_{t_1,\ldots,t_m}
\left| \prod_{i=1}^m
\Pr \left[ \frac{ \tT^{\rm in}_i
- \Exp[
\tT_i^{\rm in}]}
{(\Var[\tT_i])^{1/2}}
\leq t_i \right] - \prod_{i=1}^m \Phi(t_i)
\right| \nonumber\\ 
& \leq & C (\log \lambda)^{3d} \lambda^{-1/2} 
+ C \lambda^{-1/(2d)} (\log \lambda)^{(2d+1)/2}.
\eea
We now move on to the fifth term on the right-hand side of (\ref{0321a}).
We have 
\[ |Z_i - K_i | =   (\Var [\tilde T_i])^{-1/2} |
( T_i -\Exp [ T_i]) - ( \tT_i - \Exp[ \tT_i])| \leq (\Var[ \tilde
T_i])^{-1/2} \left( |T_i - \tilde T_i | + \Exp [|T_i - \tilde T_i |]
\right), \] 
and from just below (4.19) in \cite{py4}, we have that 
this
is bounded by $C\lambda^{-3}$ except possibly on the
set $E^c_\lambda$ which has probability less than
$C\lambda^{-3}$. Thus by (\ref{0321a}) with 
$\delta = C\lambda^{-3}$, and using (\ref{0917b})
for the second term on the right-hand side of (\ref{0321a}), we obtain 
\bea
\label{0517c}
\sup_{t_1,\ldots,t_m}
\left| \Pr \left[ \bigcap_{i=1}^m \left\{Z_i \leq t_i \right\}
\right] - \prod_{i=1}^m \Phi(t_i) \right|
 \leq  C \lambda^{-1/(2d+\eps)} 
+
C (\log \lambda)^{3d} \lambda^{-1/2} \nonumber\\
 + C \lambda^{-1/(2d)} (\log \lambda)^{(2d+1)/2}
+ C\lambda^{-3} = O (\lambda^{-1/(2d+\eps)}).
\eea
 By the
triangle inequality we have 
\bea
\label{0517d}
 \sup_{t_1,\ldots,t_m} \left|
\Pr \left[ \bigcap_{i=1}^m \left\{ \frac{T_i - \Exp[ T_i]} {( \Var [T_i])^{1/2}} \leq t_i
\right\} \right] -
\prod_{i=1}^m \Phi(t_i) \right|
\nonumber\\ \leq \sup_{t_1,\ldots,t_m} \left| \Pr \left[ \bigcap_{i=1}^m \left\{
\frac{T_i - \Exp [T_i]}
{( \Var [\tT_i])^{1/2}} 
\leq 
t_i \cdot  \left( \frac{ \Var [
T_i]}{\Var [\tT_i]} \right)^{1/2} \right\} \right] - 
\prod_{i=1}^m \Phi \left( t_i \cdot  \left( \frac{ \Var [ T_i]}{\Var
[\tilde T_i]} \right)^{1/2} \right) 
\right| \nonumber\\
+ \sup_{t_1,\ldots,t_m} \left| \prod_{i=1}^m \Phi \left(  t_i \cdot  \left( \frac{ \Var [
T_i]}{\Var [\tilde T_i]} \right)^{1/2} \right) - 
\prod_{i=1}^m \Phi(t_i) \right|. \eea
 Now from (\ref{0318c}) and (\ref{0914e}),
there is a constant $C\in(0,\infty)$
such that for all $\lambda \geq 1$ and all $t_i \in \R$
\bean \left| t_i \cdot  \left(
\frac{ \Var [ T_i]}{\Var [\tilde T_i]} \right)^{1/2} -t_i
\right| = | t_i | \left| \left( 1 + \frac{ \Var[ T_i] - \Var [\tT_i]}
{\Var [\tT_i]} \right)^{1/2} - 1 \right| \\
= |t_i | \left| \left( 1 + O( \lambda^{-3} )\right)^{1/2} - 1 \right| \leq
C |t_i| \lambda^{-3}; \eean
then since for all $s \leq t$ we
have $| \Phi(s) -\Phi(t)| \leq (t-s) \max_{s \leq u \leq t}
\varphi(u)$,
we get
\bea
\label{0917v}
 \sup_{t_i}   \left| \Phi \left(  t_i \cdot  \left(
\frac{ \Var [ T_i]}{\Var [\tilde T_i]} \right)^{1/2}
\right) - \Phi(t_i) \right| \leq C \sup_{t_i} \left(
|t_i| \lambda^{-3} \sup_{u:|u-t_i| \leq t_i C\lambda^{-3}} \varphi(u) \right)
\leq C \lambda^{-3} .
\eea
Then, considering
the second term on the right-hand side of (\ref{0517d}), we have
\bea
\label{0917w}
&  & 
\sup_{t_1,\ldots,t_m} \left| \prod_{i=1}^m \Phi \left(  t_i \cdot  \left( \frac{ \Var [
T_i]}{\Var [\tilde T_i]} \right)^{1/2} \right) - 
\prod_{i=1}^m \Phi(t_i) \right| \nonumber\\ 
& = & \sup_{t_1,\ldots,t_m} \left| \sum_{i=1}^m
\left[ \Phi \left(  t_i \cdot  \left( \frac{ \Var [
T_i]}{\Var [\tilde T_i]} \right)^{1/2} \right)
- \Phi(t_i) \right] \prod_{j=i+1}^m  \Phi(t_j) \prod_{k=1}^{i-1}
\Phi \left(  t_k \cdot  \left( \frac{ \Var [
T_k]}{\Var [\tilde T_k]} \right)^{1/2} \right) \right| \nonumber\\
& \leq & \sup_{t_1,\ldots,t_m} \sum_{i=1}^m
\left| \Phi \left(  t_i \cdot  \left( \frac{ \Var [
T_i]}{\Var [\tilde T_i]} \right)^{1/2} \right)
- \Phi(t_i) \right| \leq C \lambda^{-3} ,\eea
by (\ref{0917v}). Thus for any $\eps>0$, from (\ref{0517d}) and (\ref{0517c})
with  
(\ref{0917w}),
\bean
 \sup_{t_1,\ldots,t_m} \left| \Pr \left[ \bigcap_{i=1}^m \left\{
\frac{T_i - \Exp [T_i]}
{(\Var [T_i])^{1/2}} \leq t_i \right\}\right] - \prod_{i=1}^m \Phi(t_i) \right|
 \leq 
C\lambda^{-1/(2d+\varepsilon)} 
+ C \lambda^{-3}   
= O(  \lambda^{-1/(2d+\eps)} ).
\eean
  This completes the proof of Theorem \ref{thm1}. $\square$

\section{Indication of applications}
\label{appl}
\allco

In applying Theorem \ref{thm1}, one needs to check that the stabilization and moments 
conditions given in Definitions \ref{expstab} and \ref{momp} hold.
These conditions, or related versions thereof, 
are known to hold for many problems of interest in geometric
probability; see \cite{by2} and \cite{py4} for an indication
of problems for which exponential stabilization and moment bounds are satisfied. 

One also needs to verify the variance bound (A1): we discuss methods of doing this in Section
\ref{vars} below. In many cases, (A1) (or related versions thereof)
has been demonstrated,
see for example \cite{avbert, by2, py1}.

In Section \ref{nng} we give an example of our result as applied to the $k$-nearest neighbour graph. In particular,
we give a multivariate CLT with explicit variance scalings in the case of the nearest-neighbour
(directed) graph on disjoint subsets of the real line (Theorem \ref{nngth} below).

\subsection{Control of variances}
\label{vars}
\allco

In this section we discuss
conditions under which one can say something about the variances $\Var[T_i]$. Recall that Theorem \ref{thm1}
is stated under assumption (A1).
First we give a sufficient condition for (A1) to hold, similar in spirit to that
used by Avram and Bertsimas \cite{avbert}. 
Once again, for notational convenience we consider only the unmarked case 
with $\MM =\{1\}$.

First we introduce some notation.
Recall that  $Q_x$ denotes
the unit $d$-cube centred at $x \in \R^d$.
For a non-null Borel subset $B$ of $A \subset \R^d$ and $\lambda>0$
 we define the
following
packing of $\lambda^{1/d} B$ by unit $d$-cubes. For $\lambda^{1/d} B \subset \R^d$ let 
\bea
\label{0001}
\W_\lambda (B) := \left\{ w \in \Z^d: 
 Q_w \subseteq \lambda^{1/d} B 
 \right\},
\eea
and set
$m_\lambda(B) := \card {({\cal W}_\lambda
(B))}$. Then we define the packing ${\cal K}_\lambda (B)$ 
of $\lambda^{1/d} B$ by
 \bea
\label{0002}
{\cal K}_\lambda (B) := 
 \{ Q_w  : w \in \W_\lambda (B) \}
.
\eea

Let $f \in {\cal B}(B)$. For $w \in \W_\lambda (B)$ set
\bea
\label{Fdef}
 F_w:=F_\lambda(Q_w;B):=\sum_{x \in \Po_\lambda \cap \lambda^{-1/d} Q_w}
\xi_\lambda (x ; \Po_\lambda) \cdot f(x) .\eea
Let ${\cal F}_\lambda$ denote the $\sigma$-field generated by the points
of $\Po_\lambda$.

\begin{defn}
\label{nondeg}
Let $\{ A_w : w \in \W_\lambda(B) \}$ be a set of $m_\lambda(B)$
events in ${\cal F}_\lambda$, associated
with the $m_\lambda(B)$ cubes $Q_w(B)$, so that each $A_w$
occurs with probability uniformly bounded
away from zero. Let $J=\{w_1,\ldots,w_M\}$
be the (random) set of indices $w \in \W_\lambda (B)$ such that
$A_w$ occurs. 
Let ${\cal G}_\lambda$ denote
the $\sigma$-field generated by the random set
 $J \subseteq \W_\lambda(B)$ and the
  values of $F_\lambda(Q_w; B)$ for
$w \notin J$. 

We say that $\mulx$ is {\em nondegenerate} on $(B,f)$ if 
there exist events $\{ A_w : w \in \W_\lambda(B)\}$ 
in $\FF_\lambda$, with $\Pr (A_w) \geq \rho >0$ for all $w$,
such that: 
\begin{itemize}
\item[(i)] given ${\cal G}_\lambda$, for all $w \in J$, $\Var [ F_w | {\cal G}_\lambda ] > \eta > 0$ a.s.;
\item[(ii)] given ${\cal G}_\lambda$, for all $w,v \in J$ with $w \neq v$, $F_w$ and $F_v$
are (conditionally) independent.
\end{itemize}
\end{defn}

The idea of this condition is that the events $A_w$ essentially `isolate' cubes $Q_w$,
while allowing strictly positive variability (of the integrated measure)
within the cube, and a positive fraction
of all the cubes $Q_w$ will be so `isolated'. 

This nondegeneracy condition can often be demonstrated. In many cases, event $A_w$ will involve
a configuration of many points in an `annulus' just outside the cube $Q_w$, and an empty `moat' inside the cube,  that ensures
sufficient independence; 
see \cite{avbert} for such a construction (in a similar context)
for the total length of the $j$-th nearest-neighbour, Voronoi, and Delaunay graphs.

We now show that given the nondegeneracy condition of Definition
 \ref{nondeg}, we have lower bounds of order $\lambda$ on the variances of $T_i$.

\begin{lemm} \label{0305aa}
Let $\Gamma$ be a non-null Borel subset of $A \subset \R^d$ such that
$| \partial_r (\Gamma) | = O(r)$ as $r \downarrow 0$. 
Then, as $\lambda \to \infty$,
\bea
\label{0901ef}
 m_\lambda (\Gamma) - |\lambda^{1/d} \Gamma| =
O \left(\lambda^{(d-1)/d} \right). 
\eea
\end{lemm}
\proof This follows in a similar way to the proof of (\ref{0901e}) given previously. $\square$

 \begin{lemm}
\label{0309ee}
Suppose that $\mulx$ is nondegenerate on $(\Gamma_i,f_i)$ (see Definition
 \ref{nondeg}).
Then there exists a constant $C_i \in (0,\infty)$ such that 
for all $\lambda$ sufficiently large
\bean
 \Var [
 T_i ] \geq C_i \lambda. \eean
\end{lemm}
\proof
From the definitions of the packing
and covering defined by (\ref{0001}), (\ref{0002})
and (\ref{cover1}), (\ref{cover2}) respectively, and the
equations (\ref{0306a}) and (\ref{Fdef}),
we have that
for Borel $\Gamma \subseteq A \subset \R^d$ and $f \in \BB(\Gamma)$
\bea
\label{987}
 \langle f, \mulx \rangle
= \sum_{w \in \W_\lambda (\Gamma)} F_\lambda (Q_w;\Gamma) + \Delta_\lambda (\Gamma),\eea
where we have set
\[ \Delta_\lambda(\Gamma) := \sum_{w \in \ZZ_\lambda ( \Gamma) \setminus  \W_\lambda (\Gamma)} 
\sum_{k=1}^{N_w} 
\xi_\lambda (X_{w,k} ; \Po_\lambda ) \cdot f (X_{w,k}) \cdot
\1_\Gamma (X_{w,k}) .\]
That is, $\Delta_\lambda (\Gamma)$ gives the contributions to $\langle f, \mulx \rangle$
from cubes that are in the covering of $\lambda^{1/d} \Gamma$ but not the packing. 

Consider $\Gamma=\Gamma_i$ with $|\Gamma_i|>0$
and $f=f_i \in \BB(\Gamma_i)$. By a similar argument to (\ref{0914n}), and (\ref{0318c}), 
we have that $\Var[\Delta_\lambda (\Gamma_i)] = o(\lambda)$ as $\lambda \to \infty$. So, by
(\ref{987}) and the Cauchy-Schwarz inequality, to prove the lemma
 it suffices to show that for all $\lambda$ sufficiently large
\[ \Var \left[\sum_{w\in \W_\lambda(\Gamma_i)}  F_\lambda (Q_w; \Gamma_i)  \right] \geq C \lambda,\]
for some $C \in (0,\infty)$.

Recall Definition  \ref{nondeg}. The proof now follows the idea of Avram and Bertsimas (see
\cite{avbert}, Proposition 5).
Let $M= \sum_{w \in \W_\lambda(\Gamma_i)} \1_{A_w}$. Recall the packing
defined by (\ref{0002}). Index the cubes $Q_w$
for which $A_w$ holds
by $J=\{w_1,\ldots,w_M\} \subseteq \W_\lambda(\Gamma_i)$. Then
\bea
\label{xcvb}
 \Exp [ M] = \sum_{w \in \W_\lambda(\Gamma_i)} \Pr (A_w) \geq \rho m_\lambda(\Gamma_i) \geq
C \lambda | \Gamma_i | \geq 
 C \lambda,\eea
 using (\ref{0901ef}) for the penultimate inequality, and the fact that
 $|\Gamma_i|>0$ for the final one.
As above, let ${\cal G}_\lambda$ denote
the $\sigma$-field generated by the random set
 $J=\{w_1,\ldots,w_M\}$ and
 the values of $F_\lambda(Q_w;\Gamma_i)$ for
$w \notin J$. Then
\bean
\Var \left[  \sum_{w \in \W_\lambda(\Gamma_i)} F_\lambda (Q_w ; \Gamma_i) 
\right] \geq \Exp \left( \Var \left[ \left. \sum_{w \in J} F_\lambda (Q_w ;\Gamma_i)
+ \sum_{w \notin J} F_\lambda (Q_w ; \Gamma_i) \right| {\cal G}_\lambda \right] \right) \\
=
\Exp \left( \Var \left[ \left. \sum_{w \in J} F_\lambda (Q_w ; \Gamma_i) \right| {\cal G}_\lambda \right] \right)
, \eean
using the fact that the sum over $w \notin J$
is ${\cal G}_\lambda$-measurable. But by condition (ii) in Definition \ref{nondeg},
the $F_\lambda(Q_w ;\Gamma_i)$ for $w \in J$
are conditionally independent (under ${\cal G}_\lambda$), so we obtain
\[ \Exp \left( \Var \left[ \left. \sum_{w \in J} F_\lambda (Q_w ;  \Gamma_i) \right| {\cal G}_\lambda \right] \right)
= \Exp \sum_{w \in J} \Var [ F_\lambda (Q_w ;  \Gamma_i) | {\cal G}_\lambda ] \geq \eta 
\Exp[M],\]
by condition (i) in Definition \ref{nondeg}. Then by (\ref{xcvb}),
the proof is complete. $\square$\\

Under certain extra conditions, it is the case that
\bea
\label{bbb}
\lambda^{-1} \Var [T_i] \to \sigma_i^2 ,\eea
for some $\sigma_i^2 \in [0,\infty)$;
see \cite{by2} and \cite{mpconv}. 
Often $\sigma_i^2$ is given explicitly as an integral; however, it
is often non-trivial to compute or to
verify that it is strictly positive.

Under additional conditions
(somewhat resembling  (i) in Definition \ref{nondeg} above)
it can be shown that $\sigma_i^2>0$. 
When (\ref{bbb}) holds
with $\sigma_i^2>0$ for all $i$,
we obviously have (A1). 
Conditions of this type
were given in \cite{py1,by2}, where a form of {\em external} stabilization is used
(which roughly speaking says that not only do Poisson points beyond the radius
of stabilization for $x$ not influence $x$, but also $x$ does not influence these points). 
The results of \cite{py1,by2} 
imply that in many cases of interest (\ref{bbb}) holds with $\sigma_i^2>0$ (given extra
conditions on $f_i$ and $\kappa$). Functionals $\xi$ for which this holds include
those associated with
the total edge length of the $k$-nearest neighbour graph, and
the total number of edges in the sphere of influence graph, plus others (see \cite{py1,by2}).
Then combining (\ref{bbb}) with external stabilization and the existence of moments
(see Section 3 of \cite{py4} for some examples) one can obtain (\ref{0901d}). 

\subsection{Example: the $k$-nearest neighbour graph}
\label{nng}

The arguments indicated above are spelled
out for the particular case of the $k$-nearest neighbour graph in Section 3.1 of \cite{py4}. 
Recall that for $k \in \N$ and a locally finite point set $\X \subset \R^d$, the $k$-nearest neighbour
(undirected) graph on $\X$ (denoted ${\rm kNG}(\X)$)
is the graph with vertex set $\X$ obtained by including $\{x,y\}$
as an edge whenever $y \in \X$ is one of the $k$ nearest neighbours of $x \in \X$, or vice versa (or both).
Let $\xi(x;\X)$ be one half the sum of the lengths
in ${\rm kNG}(\X)$ incident to $x$. Thus (for example) we have that the total length of
${\rm kNG}(\X)$ is given by
\[ \sum_{x \in \X} \xi(x;\X) .\]
Suppose $\Gamma_1,\ldots,\Gamma_n$ are disjoint convex or polyhedral regions. We give two examples
of conditions on $\{f_i\}$ and $\kappa$ which, by known results together with Theorem \ref{thm1},
yield (\ref{0901d}) for this case.

First, suppose that $\kappa$ is bounded away from
$0$ on $\cup_i \Gamma_i$. Then $\xi$ is exponentially stabilizing and has moments of all orders.
If $f_i$ is continuous on $\Gamma_i$, then (\ref{bbb}) holds with $\sigma_i^2>0$ (see \cite{py4}, Section 3.1). 
Hence Theorem \ref{thm1}
applies in this case. The conditions on $f_i$ and $\kappa$ may be relaxed (see \cite{mpconv}), but then
extra work (such as making use of the nondegeneracy
argument in the present paper)
is needed to show that $\sigma_i^2>0$.

Alternatively, suppose
that $\kappa$ is equal to a positive constant $\kappa_i$ on each $\Gamma_i$, so that $\Po_\lambda$ is a homogeneous Poisson point
process
with intensity $\lambda \kappa_i>0$ on $\Gamma_i$. Suppose that $f_i =\1_{\Gamma_i}$, the indicator of $\Gamma_i$,
for each $i$. Then by the results of \cite{py1}, we again 
have that (\ref{bbb}) holds with $\sigma_i^2>0$, and so Theorem \ref{thm1} holds. In this case,
$T_i$ is the total length of ${\rm kNG}(\Po_{\lambda} \cap \Gamma_i)$.

We conclude this section by presenting an explicit multivariate 
CLT of this type, 
derived from Theorem \ref{thm1},
for the case of the nearest-neighbour (directed) graph in one dimension. 
The nearest-neighbour (directed) graph on locally finite point set $\X$
is the graph with vertex set $\X$ obtained by including $(x,y)$ as a (directed) edge
from $x \in \X$ to $y \in \X$ when $y$ is the nearest neighbour of $x$ (arbitrarily breaking any ties).
The required moments, regularity and stabilization conditions
all follow from previous work (particularly \cite{mpconv,py1}), and the fact that the limiting
variance is non-zero follows from an explicit calculation (which we give below)
 based on the general results
of \cite{mpconv}.

For a finite
set $\X \subset (0,1)$ and a Borel set $\Gamma \subseteq (0,1)$, 
let $\L^\alpha (\X;\Gamma)$ denote the
total 
weight
of the nearest-neighbour (directed) graph
on $\X$, with $\alpha$-power weighted edges, counting
only edges originating from points of $\X \cap \Gamma$. That is, if $d(x;\X):=d_2 (x; \X \setminus \{ x\})$ denotes
the (Euclidean)
 distance from $x$ to its nearest neighbour in $\X$, take 
\bea
\label{oiu}
 \xi(x;\X) = (d(x;\X))^\alpha ,\eea
for some fixed parameter $\alpha \in (0,\infty)$. Then
\[ \L^\alpha (\X;\Gamma) = \sum_{x \in \X \cap \Gamma} \xi(x;\X).\]

For $m \in \N$, 
let $\Gamma_1,\ldots,\Gamma_m$ be  disjoint, finite, non-null
interval subsets 
of $\R$.
In particular, let $\pi_i = |\Gamma_i| \in (0,\infty)$ be the length
of the interval $\Gamma_i$.
Take $f_i =\1_{\Gamma_i}$. Let the underlying density
$\kappa$ be piecewise Borel-measurable, bounded away from $0$ and from $\infty$, on each interval $\Gamma_i$; in particular, for
each $i$ set
$\kappa(x)=\kappa_i(x)$ for $x \in \Gamma_i$, where
$\kappa_i \in \BB (\Gamma_i)$ and $\kappa_i(x) >0$ for all $x \in \Gamma_i$. Consider the unmarked case (so $\MM =\{1\}$). Then
for $\lambda>0$, $\Po_\lambda$ is a Poisson point process with intensity $\kappa_i(x) \lambda$ on each $\Gamma_i$.
Using the notation of Theorem \ref{thm1}, in this set-up we have that 
\[ T_i = \langle \1_{\Gamma_i} , \mu_\lambda^\xi \rangle = \sum_{x \in \Po_\lambda \cap \Gamma_i}
\xi_\lambda (x ; \Po_\lambda)
= \sum_{x \in \Po_\lambda \cap \Gamma_i}
\xi ( \lambda x ; \lambda \Po_\lambda),\]
the final equality by translation-invariance. By the scaling properties (`homogeneity')
of $\xi$ as given by (\ref{oiu}), we have 
\[ T_i= \sum_{x \in \Po_\lambda \cap \Gamma_i}
\xi ( \lambda x ; \lambda \Po_\lambda) = 
\lambda^{\alpha} \sum_{x \in \Po_\lambda \cap \Gamma_i}
\xi ( x ; \Po_\lambda) = \lambda^{\alpha} \L^\alpha (\Po_\lambda;\Gamma_i).
\]
All relevant stabilization, regularity and moments conditions are satisfied. 
Let $\H_1$ denote a homogeneous Poisson point process of unit 
intensity on $(0,1)$, and
let $\U_n$ denote a binomial point process consisting
of $n$ independent uniform random points on $(0,1)$. Then by Theorems 2.2 and 2.4 of \cite{mpconv},
for $\alpha >0$
\bea
\label{xxx}
\lim_{\lambda \to \infty}
\lambda^{-1} \Var[T_i] 
=
 \lim_{\lambda \to \infty} \lambda^{2\alpha-1} 
\Var[\L^\alpha (\Po_{\lambda};\Gamma_i)] \nonumber\\
= V_\alpha \int_{\Gamma_i} \kappa_i(x) \ud x + \left( \delta_\alpha \int_{\Gamma_i} \kappa_i(x) \ud x \right)^2,
\eea
where\bea
\label{eek}
 V_\alpha :=
\lim_{n \to \infty} n^{2\alpha-1} 
\Var[\L^\alpha (\U_n;(0,1))]
,\eea
and
\bea
\label{eekd} \delta_\alpha := \Exp [ d(0;\H_1)^\alpha ] + \int_\R \Exp [ d(0;\H_1 \cup \{ y\})^\alpha
- d(0; \H_1)^\alpha ] \ud y.\eea

Let $\Gamma(\cdot)$ denote the (Euler) Gamma function, and let
 $_2 F_1 (\cdot, \cdot; \cdot; \cdot)$ denote
the (Gauss) hypergeometric function (see e.g.~\cite{as}, Chapter 15).
By (\ref{eek}) and equation (24) in \cite{ong}, we have that for $\alpha >0$
\bea
\label{valpha}
 V_\alpha = (4^{-\alpha} + 2 \cdot 3^{-1-2\alpha} ) \Gamma (1+2\alpha)
- 4^{-\alpha} (3+\alpha^2) \Gamma (1+\alpha)^2 \nonumber\\
+ 8 \cdot \frac{6^{-\alpha-1} \Gamma (2+2\alpha)}{(1+\alpha)}
~_2 F_1 (-\alpha, 1+\alpha; 2+\alpha; 1/3) .\eea
We now compute $\delta_\alpha$.
By standard properties of the Poisson process, $D:=d(0;\H_1)$ is distributed as
an exponential random variable with parameter $2$. So
we have that for $\alpha >0$
\[ \Exp [ D^\alpha ] = \int_0^\infty 2 r^{\alpha } \exp (-2r) \ud r = 2^{-\alpha} \Gamma (1+\alpha) ,\]
(using Euler's Gamma integral; see e.g.~6.1.1 in \cite{as}). 
By Fubini's theorem and (\ref{eekd}) we have
\bea
\label{dalpha}
\delta_\alpha = \Exp \left[ D^\alpha - 2 \int_0^D ( D^\alpha - t^\alpha ) \ud t \right]
= \Exp [ D^\alpha +( (2/(1+\alpha)) -2) D^{1+\alpha} ] \nonumber\\
= 2^{-\alpha} \Gamma (1+\alpha) - \frac{2\alpha}{1+\alpha} 2^{-1-\alpha} \Gamma (2+\alpha)
= 2^{-\alpha} \Gamma (1+\alpha) (1-\alpha),
\eea
using the functional relation $\Gamma (x)=x^{-1} \Gamma(1+x)$ (see e.g.~6.1.15 in \cite{as}) for the final
equality. 
Of note is the fact that $\delta_1=0$, so that in the $\alpha=1$ case the constant in the limiting (scaled)
variance is the same in the Poisson and binomial cases. For $\alpha \neq 1$, $\delta_\alpha^2>0$ and the variance in the Poisson case is greater than that in the binomial case, 
as one expects (the Poisson process introduces additional randomness).

Also by Theorem 2.1 of \cite{mpconv} and equation (22) in \cite{ong}, we have that for $\alpha>0$
\[ \lambda^{\alpha -1} \Exp [ \L^\alpha (\Po_\lambda ; \Gamma_i)] \to 2^{-\alpha} \Gamma (1+\alpha) \int_{\Gamma_i}
\kappa_i (x) \ud x,\]
as $\lambda \to \infty$. Thus we have the following application of Theorem \ref{thm1}.

\begin{theo}
\label{nngth}
For $m \in \N$, 
let $\Gamma_1,\ldots,\Gamma_m$ be disjoint 
intervals in $\R$ 
with $|\Gamma_i| = \pi_i \in (0,\infty)$. Let
$\kappa(x)=\sum_{i=1}^m \kappa_i(x) \1_{\Gamma_i}(x)$ where, for each $i$, $\kappa_i \in \BB (\Gamma_i)$ and
$\kappa_i(x)>0$ for all $x \in \Gamma_i$. Suppose $\alpha \in (0,\infty)$. 
\begin{itemize}
\item[(i)]
For $1 \leq i \leq m$,
\[ \lim_{\lambda \to \infty} \lambda^{\alpha-1} \Exp [ \L^\alpha (\Po_\lambda ; \Gamma_i) ] = 
 2^{-\alpha} \Gamma (1+\alpha) \int_{\Gamma_i} \kappa_i (x) \ud x.\]
\item[(ii)]
For $1 \leq i \leq m$,
\[ \lim_{\lambda \to \infty} \lambda^{2\alpha-1} 
\Var[\L^\alpha (\Po_{\lambda};\Gamma_i)] \nonumber\\
= V_\alpha \int_{\Gamma_i} \kappa_i(x) \ud x + \left( \delta_\alpha \int_{\Gamma_i} \kappa_i(x) \ud x \right)^2
=: \sigma_i^2,\]
where $V_\alpha$ and $\delta_\alpha$ are given by (\ref{valpha}) and (\ref{dalpha}) respectively.
\item[(iii)]
Given  $\eps >0$,
there exists $C \in (0,\infty)$ such that
for all $\lambda \geq 1$,
\[ \sup_{t_1, \ldots, t_m \in \R} \left| \Pr \left[ \bigcap_{i=1}^m \left\{ \frac{\L^\alpha (\Po_\lambda ; \Gamma_i)
- \Exp [ \L^\alpha (\Po_\lambda ; \Gamma_i) ]}
{ ( \Var [ \L^\alpha (\Po_\lambda ;\Gamma_i ) ])^{1/2} } \leq t_i \right\}
\right] - \prod_{i=1}^m \Phi (t_i) \right| \leq C \lambda^{\eps-(1/2)}.\]
\end{itemize}
\end{theo}

Part (iii) of Theorem \ref{nngth} is our multivariate CLT. In
 the particular case of piecewise constant $\kappa$, that is $\kappa_i(x) = \kappa_i \in (0,\infty)$ for all $x \in \Gamma_i$, we have that
\[ \int_{\Gamma_i} \kappa_i (x) \ud x = \kappa_i | \Gamma_i | = \kappa_i \pi_i,\]
and so, for example, $\sigma_i^2 = V_\alpha \kappa_i \pi_i + \delta_\alpha^2 \kappa_i^2 \pi_i^2$.
Table \ref{tabvarpo} gives some values of the constants
$V_\alpha$, given by (\ref{valpha}), and $\delta_\alpha^2$, given by (\ref{dalpha}).
\begin{table}[!h]
\begin{center}
\begin{tabular}{|l|l|l|l|l|l|}
\hline
$\alpha$  &  1/2 & 1 & 2 & 3 & 4   \\
\hline
$V_\alpha$ &    $\frac{1}{2} + \sqrt{2} \arcsin (1/\sqrt{3})
- \frac{13\pi}{32} \approx 0.094148$ & $\frac{1}{6}$ &  $\frac{85}{108}$ & $\frac{149}{18}$ & $\frac{135793}{972}$  \\
\hline
$\delta_\alpha^2$ &  $\frac{\pi}{32}$ & 0 & $\frac{1}{4}$ & $\frac{9}{4}$ & $\frac{81}{4}$ \\
\hline
\end{tabular}
\caption{Some values of $V_\alpha$ and $\delta_\alpha^2$.}
\label{tabvarpo}
\end{center}
\end{table}

\begin{center} \textbf{Acknowledgements} \end{center}
The authors began this work while at the University of
Durham, where AW was
 supported by an EPSRC studentship. Some of this work was carried out while AW was
 at the University of Bath. MP thanks the Institute for Mathematical Sciences
 at the National University of Singapore for its hospitality.

\end{document}